\documentclass[11pt]{amsart}
\usepackage{comment,amsmath,amssymb}
\usepackage[enableskew]{youngtab}
\usepackage{young}
\title{Equality of Schur's Q-functions and their skew analogues}
\author{Hadi Salmasian}
\date{February 1, 2006}
\address{Department of Mathematics and Statistics\\
Queen's University\\
Jeffery Hall, University Avenue\\
Kingston, ON K7L 3N6\\ Canada
}
\email{hadi@mast.queensu.ca}

\newtheorem{example}{Example}
\newtheorem{theorem}{Theorem}
\newtheorem{definition}[theorem]{Definition}
\newtheorem{corollary}[theorem]{Corollary}
\newtheorem{proposition}[theorem]{Proposition}
\newtheorem{lemma}[theorem]{Lemma}

\newcommand{\ton}{\overline{\textrm{\small1}}}
\newcommand{\ttw}{\overline{\textrm{\scriptsize2}}}
\newcommand{\tth}{\overline{\textrm{\scriptsize3}}}
\newcommand{\tfo}{\overline{\textrm{\scriptsize4}}}
\newcommand{\tfi}{\overline{\textrm{\scriptsize5}}}
\newcommand{\tsi}{\overline{\textrm{\scriptsize6}}}
\newcommand{\tse}{\overline{\textrm{\scriptsize7}}}

\newcommand{\BB}{\mathbb{B}}

\newcommand{\DD}{D_{\lambda/\mu}}
\newcommand{\ov}{\overline}
\begin{document}

\begin{abstract}
 We find a simple criterion for the equality
$Q_\lambda=Q_{\mu/\nu}$ where $Q_\lambda$ and
$Q_{\mu/\nu}$ are Schur's Q-functions on infinitely
many variables.
\end{abstract}

\maketitle

\section{Introduction}

Schur's Q-functions are very interesting
analogues of the (standard)
Schur functions $s_\lambda$
in several combinatorial and
representation-theoretic contexts. Examples
of their analogy include the shifted RSK
correspondence, the shifted Littlewood-Richardson
rule, and the character theory of representations
of queer Lie superalgebras. In this note we
study when certain shifted Littlewood-Richardson
coefficients (the $f_{\mu,\nu}^\lambda$ in the
language of \cite{st1}) are zero or one. 
Studying questions of the same nature has been of interest to a 
number of authors. In particular, one should mention Stembridge's
recent classification of multiplicity-free products of Schur functions
\cite{st2} which was generalized to P-functions in \cite{wi1}. 
The related
questions as to when two ribbon Schur functions are
equal and when a Schur function is equal to a skew
Schur function were answered in \cite{wi2}
and \cite{wi3}. Here we show
that the latter problem has a simple answer
for Schur's Q-functions as well.\\

This note is organized as follows. In the next section
we give all the
required definitions. In the third section we prove our
main result.\\

{\noindent\bf Acknowledgement.}
The author thanks the referee for reading the paper very carefully.

\section{Schur's Q-functions and shifted tableaux}
Our notation is compatible with Stembridge's paper \cite{st1}.
A strictly decreasing sequence $\lambda=\{\lambda_1>...>\lambda_k\}$
of positive integers is called a {\it distinct partition} of $n$
if the sum of
the $\lambda_i$'s is equal to $n$.
The $\lambda_i$'s are called the {\it parts} of $\lambda$.
A partition is represented
by a shifted Young diagram as follows:
there
are  $\lambda_i$ boxes in the $i$-th row, after $i-1$ empty positions.
We denote this diagram by $D_\lambda$.
\begin{example}\label{one}
\rm
Let $\lambda=\{7>4>2>1\}$. Then the diagram $D_\lambda$ which represents
$\lambda$ is given below.
$$\young(\ \ \ \ \ \ \ ,:\ \ \ \ ,::\ \ ,:::\ )\vspace{-4mm}$$
\end{example}
\vspace{-3mm}
Consider an ordered alphabet
$$\mathcal{A}=\{\overline{1}<1<\overline{2}<2<\overline{3}<3<\cdots\}.$$
The letters $\overline{1},\overline{2},...$
will be referred to as {\it marked}, whereas the letters
$1,2,...$
will be referred to as {\it unmarked}.

\begin{definition}\label{GSYT}
By a generalized shifted Young tableau (GSYT) of
shape $D_\lambda$ we mean a
filling of a given Young diagram $D_\lambda$ with
letters from $\mathcal A$
such that
the following properties hold:

\begin{enumerate}
\item[$\bullet$] The rows and columns are weakly increasing.
\item[$\bullet$] Each row contains each marked letter at most once.
\item[$\bullet$] Each column contains each unmarked letter
at most once.
\end{enumerate}
\end{definition}
Let $T$ be a given GSYT. We define $x^T$ to be the monomial
$x_1^{a_1}x_2^{a_2}\cdots$ where $a_s$ is the total number of 
occurrences
of $s$ or $\overline{s}$ in $T$. Schur's Q-function $Q_\lambda$ is equal
to
\begin{equation}\label{qschur}
\sum_{T}x^T
\end{equation}
where the sum is over all GSYT of shape $D_\lambda$. The sequence 
$a_1,a_2,...$ is called the {\it content} of $T$.\\

Now suppose $\lambda$ and $\mu$ are two distinct partitions. Assume that
for each $i$, $\lambda_i\geq \mu_i$, so that $D_\mu$ lies inside $D_\lambda$.
The {\it shifted skew diagram} associated
to $\lambda/\mu$ is the set-theoretic difference of
$D_\lambda$ and $D_\mu$. It is represented by $D_{\lambda/\mu}$.

\begin{example}\label{two}\rm
Let $\lambda$ be as in Example \ref{one} and let $\mu=\{4>2\}$.
Then $D_{\lambda/\mu}$ is represented by the following skew Young
diagram:
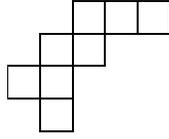
\vspace{-2mm}
$$\young(::::\ \ \ ,:::\ \ ,::\ \ ,:::\ )$$
\end{example}
\vspace{-3mm}
We may assume that $\DD$ is placed on the Cartesian plane in the usual way
such
that the centers of boxes lie on the lattice
of points with integer coordinates. Let $\BB_{x,y}$ denote the box
whose center is at the point with coordinates $(x,y)$. The following 
easy lemma includes some basic properties of shifted skew diagrams.

\begin{lemma}\label{propD}
Let $\DD$ be an arbitrary shifted skew diagram on the Cartesian plane.
For any integer $y$ define 
\vspace{-3mm}
$$
R_y=\{\BB_{x,y}|\BB_{x,y}\in\DD\}.
$$  
\begin{enumerate}
\item[$\bullet$]
If $\BB_{u,v}\in\DD$ and $R_{v+1}\neq\emptyset$ then there exists
an integer $t\geq u$ such that $\BB_{t,v+1}\in R_{v+1}$.
\item[$\bullet$] 
For a fixed $v$ such that $R_v\neq\emptyset$, 
let $u$ be the smallest number for which 
$\BB_{u,v}\in R_v$. Assume $\BB_{u-1,v+1}\in\DD$. Then for any 
$v'\leq v+1$ such that $R_{v'}\neq\emptyset$, the 
following statement is
true:
\begin{itemize}
\item[$\diamond$] $\BB_{u+v-v',v'}\in R_{v'}$ and for any $s$ if
$\BB_{s,v'}\in R_{v'}$ then $s\geq u+v-v'$.
\end{itemize}
\end{enumerate}
\end{lemma}
Schur's skew Q-function $Q_{\lambda/\mu}$ is equal to a summation
similar to (\ref{qschur}), where the summation is now
on all shifted skew tableaux $T$, with underlying diagram
$\DD$
filled by the alphabet $\mathcal A$, such that they satisfy
the properties of Definition \ref{GSYT}. The function 
$Q_{\lambda/\mu}$
can be expressed as a linear combination
of the functions $Q_\nu$ for various
$\nu$ as follows. We have
\vspace{-2mm}
\begin{equation}\label{skewQR}
Q_{\lambda/\mu}=\sum_{\nu} f_{\mu\nu}^\lambda Q_\nu
\vspace{-3mm}
\end{equation}
where the summation is over all distinct partitions $\nu$.
Here $f_{\mu\nu}^\lambda$ is the number of {\it amenable} tableaux
of shape $D_{\lambda/\mu}$ and content $\nu$. 
We define the amenable tableaux in
Definition \ref{amenable} below.
However, before doing so, we need some notation. For a given
(possibly skew) GSYT such as $T$, the row word of $T$ is the word
obtained by reading the rows of $T$ consecutively
from left to right starting
with the bottom row. We denote the row word of $T$ by 
$w(T)$.
Now, let $w=w_1\cdots w_n$ be an arbitrary word of length $n$ such that
for any $s$ we have
$w_s\in\mathcal A$.
Define a function $m_i(j)$ as follows.
$$
m_i(j)=\left\{\begin{array}{ll}
\textrm{number of times } i\textrm{ appears}&\textrm{    if }1\leq j\leq n\\
\textrm{among } w_{n-j+1},...,w_n \\
\ & \ \\
m_i(n)+\textrm{number of times }\overline{i}\textrm{ appears} &
\textrm{    if }n+1\leq j\leq 2n\\
\textrm{among }
w_1,...,w_{j-n}\\
\end{array}\right.\vspace{-1.5mm}$$
By convention, we assume $m_i(0)=0$ for any $i>0$.
\vspace{-1.5mm}
\begin{definition}\label{amenable}
Let $k>1$ be an integer.
A word $w=w_1\cdots w_n$ is called $k$-amenable iff it satisfies 
the following properties:
\begin{enumerate} 
\item[$\bullet$] For any $j\in\{0,...,n-1\}$ 
if $m_k(j)=m_{k-1}(j)$
then $w_{n-j}\notin\{k,\overline{k}\}$.
\item[$\bullet$] For any $j\in\{n,...,2n-1\}$, 
if $m_k(j)=m_{k-1}(j)$
then $w_{j-n+1}\notin\{\overline{k},k-1\}$.
\item[$\bullet$] If $j$ is the smallest number such that
$w_j\in\{k,\ov{k}\}$, then $w_j=k$.
\item[$\bullet$] If $j$ is the smallest number such that
$w_j\in\{k-1,\ov{k-1}\}$, then $w_j=k-1$.
\end{enumerate}

A word $w$ is called amenable if it is $k$-amenable for any $k>1$.

\end{definition}

\noindent{\bf Remark.} Suppose $w=w_1\cdots w_n$ is $k$-amenable for some 
$k>1$. 
Then it follows from 
Definition \ref{amenable} that 
if 
$m_{k-1}(2n)>0$ then $m_k(2n)<m_{k-1}(2n)$. (To prove this, 
first we show that Definition \ref{amenable} implies 
$m_{k-1}(j)\geq m_k(j)$ for any $j\in\{1,...,2n\}$. Then 
we pick the largest $j$ such that $w_j\in\{k,\overline{k}\}$,
and we show that Definition \ref{amenable} implies that
there must exist a $j_1>j$ such that $w_{j_1}=\overline{k-1}$. 
The details of the argument are 
left to the reader.) Consequently, if $w=w_1\cdots w_n$ is amenable, then
$$
m_1(2n)\geq m_2(2n)\geq m_3(2n)\geq \cdots
$$

\begin{definition}
A given GSYT is called amenable iff its row word is amenable.
\end{definition}

\section{The main result}
In this section we prove the main statement of this note.

\begin{definition}
A shifted skew diagram $\DD$ is called strange iff
$Q_{\lambda/\mu}=Q_\nu$ for some distinct partition $\nu$.
\end{definition}

\begin{theorem}\label{main}
$\DD$ is a strange diagram if and only if $\lambda/\mu=\overline{\lambda}/\overline{\mu}$
where 
\begin{itemize}
\item[$\bullet$] $\overline{\lambda}$ is arbitrary and $\overline{\mu}=\{\}$.
\item[$\bullet$] 
$\overline{\lambda}=\{m>m-1>\cdots>1\}$ and 
$\overline{\mu}=\{\mu_1>\cdots>\mu_l\}
\textrm{ where } 0<l<m-1.$
\item[$\bullet$] $\overline{\lambda}=\{p+q+r>p+q+r-1>p+q+r-2>\cdots>p\}$
and \\
$\overline{\mu}=\{q>q-1>\cdots>1\}$
where $p,q,r$ are integers such that $p,q\geq 1$, $r\geq 0$. 
\item[$\bullet$]
$\overline{\lambda}=\{p+q>p+q-1>\cdots>p+q-r\}$ and $\overline{\mu}=\{q>q-1>\cdots>q-r\}$
where $p,q,r$ are integers such that $p>0$ and $q>r\geq 0$.

\end{itemize}

\end{theorem}

\noindent{\bf Remark.} The reader 
should note that there are overlaps among the
cases for special values of $p,q,r$. Moreover, whether or not 
$\DD$ is strange 
only depends on $\lambda/\mu$ and not on $\lambda$ and $\mu$ individually. 
Theorem \ref{main} identifies $\DD$ by identifying
all possible differences $\lambda/\mu$, but not all possibilities 
of $\lambda$ and $\mu$.
The latter problem is not hard once we have Theorem \ref{main}.
\vspace{-2mm}\\

\vspace{-1.5mm}
The proof of Theorem \ref{main} will be given throughout this section.
From 
equation (\ref{skewQR}) it follows that $\DD$ is strange 
if and only if there exists
a unique amenable filling of $D_{\lambda/\mu}$. Our approach is
to rule out various possibilities for the shape of $D_{\lambda/\mu}$
by demonstrating the existence
of at least two different amenable fillings in each case.\\

\subsection{An algorithm for finding an amenable filling}
\label{algorithm}

We give a simple algorithm to construct an amenable tableau 
of
any given 
shape $D_{\lambda/\mu}$. 
The output of the algorithm
is an amenable tableau of content $\nu$ for some distinct partition 
$\nu$. Note that $\nu$ is generated by the algorithm and is not an
input.\\

\noindent {\bf Notation.} Let $\BB_{x,y}$ be a box in a diagram. 
The operation of filling (or replacing the current entry of) 
$\BB_{x,y}$ by the letter $w\in\mathcal{A}$
is represented by 
$\BB_{x,y}\leftarrow w$. The operation of filling (or replacing 
the entry of) $\BB_{x,y}$ by the current entry of $\BB_{x',y'}$
is represented by $\BB_{x,y}\leftarrow\BB_{x',y'}$.\\

The algorithm is given below.\\

\begin{enumerate}
\item[\tt Step 1.] Set 
$k=1$ and $\DD^{(1)}=\DD$.

\item[\tt Step 2.] Set
$P_k=\{\BB_{x,y}\,| \BB_{x,y}\in\DD^{(k)} \textrm{ but }
\BB_{x-1,y+1}\notin \DD^{(k)}\}$.
\item[\tt Step 3.] Put a $k$ or a $\ov{k}$ in any of the boxes
in $P_k$ according to the following rule:
$$
\left\{
\begin{array}{l}
\BB_{x,y}\leftarrow k\textrm{ if }\BB_{x,y-1}\notin P_k\\
\BB_{x,y}\leftarrow \ov{k}\textrm{ if }\BB_{x,y-1}\in P_k
\end{array}
\right.
$$
\item[\tt Step 4.] Remove all boxes of $P_k$ from $\DD^{(k)}$.
Let $\DD^{(k+1)}$ be the diagram obtained after removing
boxes. If $\DD^{(k+1)}$ has no boxes, then stop.
\item[\tt Step 5.] Increase $k$ by one. Go back to Step 2.\\

\end{enumerate}

Let $k$ be a positive integer. 
Then each connected component of $P_k$ 
forms a ``path'' of boxes in a connected component of the diagram
$\DD^{(k)}$.
(A connected component 
of $\DD^{(k)}$ (or $P_k$)
is a 
maximal subset of
boxes of $\DD^{(k)}$ (or $P_k$) which can be
ordered in a sequence so that each box has a 
common edge with at least
one of the boxes preceding it in the sequence.)
Each of these paths can be directed 
as follows. Let $Q$ be a connected component of $P_k$.
We know that $Q$ is a path of boxes. The 
first box of this path is the box $\BB_{x,y}\in Q$ such that
$\BB_{x,y+1}\notin Q$ and $\BB_{x+1,y}\notin Q$. The last
box of $Q$ is the box $\BB_{x,y}\in Q$ such that 
$\BB_{x-1,y}\notin Q$ and $\BB_{x,y-1}\notin Q$.

\begin{example}\label{exexloc}\rm
The first and last boxes of the following path are marked with a 
dot and a cross respectively. It is traversed from the dotted box
to the crossed one.
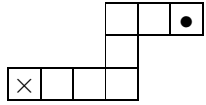
\vspace{-1mm}
$$
\young(:::\ \ \bullet,:::\ ,\times\ \ \ )
\vspace{-2.5mm}
$$
\end{example}

\begin{lemma}\label{paths}
For any $k>1$, $P_k=\{\BB_{x,y}|\BB_{x,y}\in\DD\textrm{ and }
\BB_{x-1,y+1}\in P_{k-1}\}$.
\end{lemma}
\vspace{-2mm}

\begin{example}\rm
Let $\lambda=\{7>5>3>2>1\}$ and $\mu=\{4>1\}$. Then the algorithm
provides the following amenable filling.
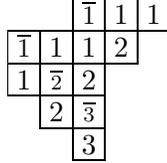
\vspace{-3mm}
$$
\begin{array}{c}
\young(::::\ton11,::\ton112,::1\ttw2,:::2\tth,::::3)
\end{array}
\vspace{-5mm}
$$
\end{example}

\begin{lemma}\label{propT}
Let $k>1$ and $\BB_{x,y}\in P_k$. Then
\begin{itemize}
\item[$\bullet$] $\BB_{x,y+1}\in\DD$ and
$\BB_{x,y+1}\in P_k\cup P_{k-1}$. 
\item[$\bullet$] If 
$\BB_{x,y+1}\in P_{k-1}$ then $\BB_{x-1,y+1}\in P_{k-1}$ as well.
\end{itemize}

\end{lemma}

\begin{proof}
Left to the reader!
\end{proof}

The following lemma is a very simple but useful criterion of
amenability. Its proof is left to the reader.

\begin{lemma}\label{checktwo}
Let $T$ be a GSYT of shape $\DD$ and let $w=w(T)$ be the row
word of $T$. For any integer $k>1$, form the word
$w^{(k)}$ from $w$ by dropping all the letters in $w$
which are not in $\{k,\overline{k},k-1,\overline{k-1}\}$.
(For example, if $w=2\ov{2}3\ov{2}211\ov{2}\,\ov{1}11$ then
$w^{(3)}=2\overline{2}3\overline{2}2\overline{2}$.)
Then $w$ is amenable if and only if $w^{(k)}$ is either $k$-amenable
or empty for any $k>1$.

\end{lemma}

\begin{lemma}\label{amenamen}
The filling of the boxes obtained by the previous algorithm
is amenable.
\end{lemma}

\begin{proof}

Properties of Definition \ref{GSYT} are
satisfied trivially. It remains to show that
this filling is amenable. By Lemma \ref{checktwo}
it suffices to check $k$-amenability
of $w(T)^{(k)}$ for any $k$.
It is easily seen that if
$\DD$ is disconnected but each connected component
of $T$ is amenable, then $T$ is amenable too. 
Therefore we can assume $\DD$ is connected.
Next we show that
$w(T)^{(2)}$ is amenable. The argument in the general case is
similar. Let the length of $w(T)^{(2)}$ be $n$.
Amenability of $w(T)^{(2)}$ follows from the
following facts:
\begin{enumerate}
\item[$\bullet$]  For any box filled by a $2$,
there is a box above it in the same column
which is filled by a $1$. This implies that
$m_1(j)\geq m_2(j)$ for any
$j\leq n$.
\item[$\bullet$]
Since the last entry of the path $P_1$ is filled by a $1$ and
there is no $2$ below it, we have $m_1(n)>m_2(n)$.
\item[$\bullet$] There is exactly one $\ov{1}$ in every row of $\DD$
in which $P_1$ has boxes,
except
for the lowest row among them.
The same statement holds for $P_2$ and $\ov{2}$.
\item[$\bullet$] 
Set $X_i=\{y |\textrm{ for some }x, \BB_{x,y}\in P_i\}$. Then 
$$
X_1=\{s\,| s\textrm{ is an integer and } c_1\leq s\leq c_2\}
$$ 
for integers $c_1\geq c_2$. Moreover
$X_2\subseteq\{s-1\, | s\in X_1\}$. 
\end{enumerate}
\end{proof}
\subsection{Disconnected diagrams}
\label{disconnected}
A slight modification of the algorithm of section \ref{algorithm}
can be applied to show that $\DD$ is not a strange diagram
whenever 
there exist disconnected $P_k$'s.
More accurately, we have the following lemma.

\begin{lemma}\label{discon}
Let $T$ be a GSYT obtained by applying the algorithm of section
\ref{algorithm} to $\DD$. Suppose that for some $k$, $P_k$ is 
disconnected.
Then $\DD$ has an amenable filling different from the one
given by the algorithm.
\end{lemma}
\begin{proof}
Let $P_k=P^1\cup\dots\cup P^l$, $l>1$, where $P^i$'s are 
mutually disjoint
paths, ordered such that for any $i$,
$P^i$ lies to the northeast of $P^{i+1}$. Note that the $P^i$'s
belong to mutually disjoint connected components of 
the diagram $\DD^{(k)}$.
Since the filling of $P_k$ is
obtained by the algorithm, the last box of $P^1$ is filled by
a $k$. We claim that
if we change it to a $\ov{k}$, the new filling is still amenable.
To this end we use Lemma
\ref{checktwo}.
One can see that changing a $k$ to a $\ov{k}$ does not affect 
$k$-amenability
of $w(T)^{(k)}$ unless we are changing the last box of $P^l$.
Moreover $w(T)^{(k+1)}$ remains $k+1$-amenable
since
$w(T)^{(k+1)}$ can be written as $w(T)^{(k+1)}=w^1\cdots w^l$
where $w^i$ is the row word of the intersection
of $P_k\cup P_{k+1}$ with the connected component 
of $\DD^{(k)}$ which contains $P^{l-i+1}$. 
A bookkeeping argument using 
the facts that the words
$w^1,...,w^{l-1}$ are
$k+1$-amenable and $w^l$ is ``not far from'' 
being $k+1$-amenable
completes the proof. 
\end{proof}

From Lemma \ref{discon} we conclude that if $\DD$ 
is a strange diagram 
then for any $k$ the set
$P_k$ which is obtained by the algorithm of section
\ref{algorithm}
is in fact a (connected)
path. This
fact simplifies the case-by-case analysis of possible shapes of $\DD$.
Probably a few remarks are necessary before starting the next 
section:

\vspace{1mm}
{\noindent$\bullet$\bf Amenability arguments.}
Throughout the next section, at several points we 
give procedures to
modify an existing filling and
claim that the new filling will be amenable too.
The proofs of amenability of these new fillings
are very similar in nature.
However, most of these amenability proofs are 
not given because they are tedious but very simple.
The main ideas are using Lemma \ref{checktwo} and 
bookkeeping arguments. 

\vspace{1mm}
{\noindent $\bullet$\bf Figures.} 
Throughout the next section there are several 
illustrations which help the reader understand 
the effect of
procedures on the fillings. The figures
usually demonstrate a union of 
$P_m,P_{m-1},...,P_r$ for some $r$. Only the
boxes whose entries change by
the procedures are shown. The dots show that there
may or may not be boxes in the direction of the dots.
See the following example.

\begin{example}\rm
Let $m$ denote the largest integer such that
$P_m\neq\emptyset$.
Suppose $m=7$ and the union $P_6\cup P_7$
is illustrated by the diagram
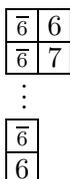
\vspace{-2mm}
$$ 
\begin{array}{ll}
\young(\tsi6,\tsi7)\vspace{-2mm}\\
\,\,\,\vdots\\
 \young(\tsi,6)
\end{array}
\vspace{-3mm}
$$
This
means that there may or may not be a 
vertical path of boxes 
which belong to $P_6$ in the location of the dots.
For instance $P_6\cup P_7$
may actually be
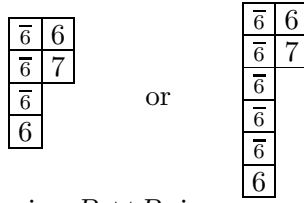
\vspace{-1mm}
$$ 
\begin{array}{ll}
\young(\tsi6,\tsi7,\tsi,6)\vspace{-2mm}\\
\,\,\,\\
\end{array}
\qquad\textrm{or}\qquad 
\begin{array}{ll}
\young(\tsi6,\tsi7,\tsi,\tsi,\tsi,6)\vspace{-2mm}\\
\end{array}
\vspace{-1mm}
$$
\vspace{-2mm}
Now again suppose $m=7$ and
the union $P_6\cup P_7$ is represented
by 
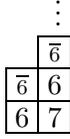
\vspace{-1mm}
$$ 
\begin{array}{l}
\ \ \ \ \ \vdots\\
\young(:\tsi,\tsi6,67)
\end{array}
\vspace{-2mm}
$$
This
means that if the path $P_6$ extends from the top right corner
of $P_6\cup P_7$,
then its direction will be vertical in the beginning. (However, 
its direction does not have to remain vertical
all the time.) For 
instance, $P_6\cup P_7$ may actually be 
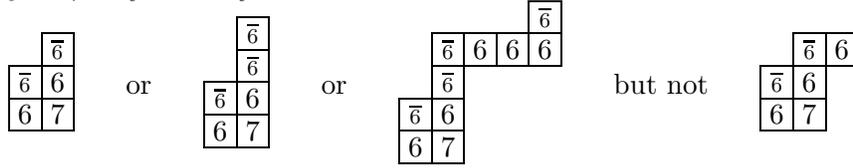
\vspace{-2mm}
$$ 
\begin{array}{l}
\young(:\tsi,\tsi6,67)
\end{array}
\quad
\textrm{  or  } 
\quad
\begin{array}{l}
\young(:\tsi,:\tsi,\tsi6,67)
\end{array}
\quad
\textrm{  or  }
\quad
\begin{array}{l}
\young(::::\tsi,:\tsi666,:\tsi,\tsi6,67)
\end{array}
\quad
\textrm{  but not  }
\quad
\begin{array}{l}
\young(:\tsi6,\tsi6,67)
\end{array}
$$
\vspace{-3mm}
Again let $m=7$ and suppose $P_6\cup P_7$
is illustrated by 
the diagram 
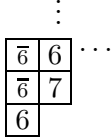
\vspace{-.5mm}
$$\begin{array}{ll}
\ \ \ \ \ \ \vdots\vspace{-1mm}\\
&\!\!\!\!\ldots\vspace{-.3cm}\\
\ \young(\tsi6,\tsi7,6)
\end{array}
\vspace{-2mm}
$$
This means that $P_6$ may extend either horizontally or 
vertically (but clearly not both, since $P_6$ is a path).
\end{example}

\subsection{A case-by-case examination}
\label{case-by-case}
Let $T$ be the filling of $\DD$ obtained by the algorithm of
section \ref{algorithm}. By Lemma \ref{discon}, $P_k$ is a connected 
path for
any $k$. Let $m$ be the largest integer for which $P_m$ is nonempty.
Then the following possibilities exist for $P_m$:

\begin{itemize}
\item[$\bullet$] $P_m$ consists of a single box.
\item[$\bullet$] $P_m$ lies within a single row
and has at least two boxes.
\item[$\bullet$] $P_m$ 
lies within a single column and has at least two boxes.
\item[$\bullet$] $P_m$ is none of the above.
\end{itemize}

We will study the first three possibilities in subsequent
sections. Here we show that the fourth case is actually impossible
when $\DD$ is a strange diagram.

\begin{lemma}
Suppose $\DD$ is a strange diagram. Then the boxes of
$P_m$ lie within a single row or column.
\end{lemma}
\begin{proof}
The following figures show how to obtain a new amenable 
filling of $\DD$ when
$P_m$ has a ``turning point''. If the path $P_m$ changes its direction
at some point clockwise, i.e. if for integers
$x,y$ we have 
$\BB_{x,y},\BB_{x,y+1},\BB_{x-1,y}\in P_m$
then all we need to do is 
$$
\BB_{x,y}\leftarrow m+1
\textrm{ and }\BB_{x,y+1}\leftarrow m.$$ 
If it changes its direction counterclockwise, i.e. if
$\BB_{x,y},\BB_{x,y-1},\BB_{x+1,y}\in P_m$,
then
we should find the smallest $y'$ such that $\BB_{x,y'}\in P_m$
and then do $\BB_{x,y'}\leftarrow m+1$ and 
$\BB_{x,y'+1}\leftarrow m$.

\begin{example}\rm
Suppose $m=5$. Then $P_m$ and the way it changes are
illustrated 
below.
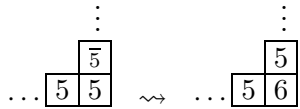
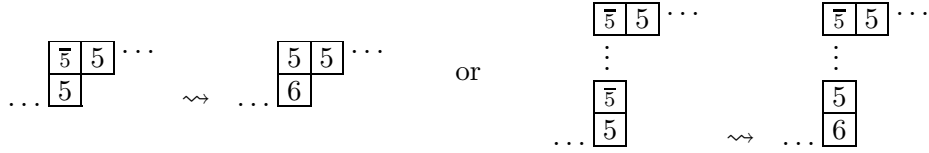
\vspace{-7mm}
$$
\begin{array}{rcr}\\
\!\!\!\!\vdots\,\,&&\!\!\!\!\vdots\,\,\\
\cdots\young(:\tfi,55)&\leadsto&
\cdots\young(:5,56)
\end{array}
$$
\begin{center}
\bf clockwise turn
\end{center}
$$
\begin{array}{rllcrll}
&&\!\!\!\!\ldots&&&&\!\!\!\!\ldots\vspace{-.3cm}\\
\ldots\!\!\!\!&\young(\tfi5,5)&&
\leadsto&
\ldots\!\!\!\!&\young(55,6)&
\end{array}
\qquad
\textrm{or}
\qquad
\begin{array}{rllcrll}
&&\!\!\!\!\ldots&&&&\!\!\!\!\ldots\vspace{-.3cm}\\
&\young(\tfi5)&&&&\young(\tfi5)\vspace{-2mm}\\
&\,\,\vdots&&&&\,\,\vdots\\
\ldots\!\!\!\!&\young(\tfi,5)&&
\leadsto&
\ldots\!\!\!\!&\young(5,6)&
\end{array}
$$

\begin{center}
\bf Counterclockwise turn
\end{center}

Checking amenability is easy and left to the reader.

\end{example}

\end{proof}

In the subsequent sections we address the other possibilities
for $P_m$.

\subsection{Case I: $P_m$ consists of one box only}
\label{CASE_I}
We will assume $\DD$ is a strange diagram and
$m>1$. Let $\BB_{x,y}$ be the box in $P_m$. Then
by Lemma \ref{propT} and Lemma \ref{paths},
$\BB_{x,y+1}$ contains $m-1$ and $\BB_{x-1,y+1}\in P_{m-1}$.
Now we have the following two possibilities:

\begin{enumerate}

\item[\bf 1.] $\BB_{x-1,y}\notin\DD$. Then it follows that $\BB_{x-1,y+1}$
should be the last box of $P_{m-1}$, as $P_{m-1}$ cannot proceed to
$B_{x-2,y+1}$ by Lemma \ref{propD}. 
This in turn yields three new possibilities:

\begin{itemize}
\item[\bf Case a.] 
$\BB_{x,y+2},\BB_{x+1,y+1}\notin P_{m-1}$. In this case $P_{m-1}$ has
two boxes only.
\item[\bf Case b.] 
$\BB_{x,y+2}\in P_{m-1}$.
\item[\bf Case c.] 
$\BB_{x+1,y+1}\in P_{m-1}$.
\end{itemize}
\item[\bf 2.] $\BB_{x-1,y}\in\DD$. Then $\BB_{x-1,y+1}$
should be filled by $\ov{m-1}$ and $\BB_{x-1,y}\in P_{m-1}$. 
One of the following six cases may take place.

\begin{itemize}
\item[\bf Case d.] $\BB_{x-2,y}\in P_{m-1}$.
\item[\bf Case e.] $\BB_{x-1,y-1},\BB_{x-1,y-2}\in P_{m-1}$.
\item[\bf Case f.] $\BB_{x-1,y-1}\in P_{m-1}$ but $\BB_{x-1,y-2}\notin\DD$.
\item[\bf Case g.] $\BB_{x-2,y},\BB_{x-1,y-1}\notin P_{m-1}$ but 
$\BB_{x+1,y+1}\in P_{m-1}$.
\item[\bf Case h.] $\BB_{x-2,y},\BB_{x-1,y-1}\notin P_{m-1}$
but $\BB_{x,y+2}\in P_{m-1}$.
\item[\bf Case i.] $\BB_{x-2,y},\BB_{x-1,y-1},\BB_{x,y+2},
\BB_{x+1,y+1}\notin P_{m-1}$.\\
\end{itemize}
\vspace{-3mm}
\end{enumerate}
\noindent Assuming $m=7$ we can illustrate $P_m\cup P_{m-1}$ 
in cases a to i by
the following figures.
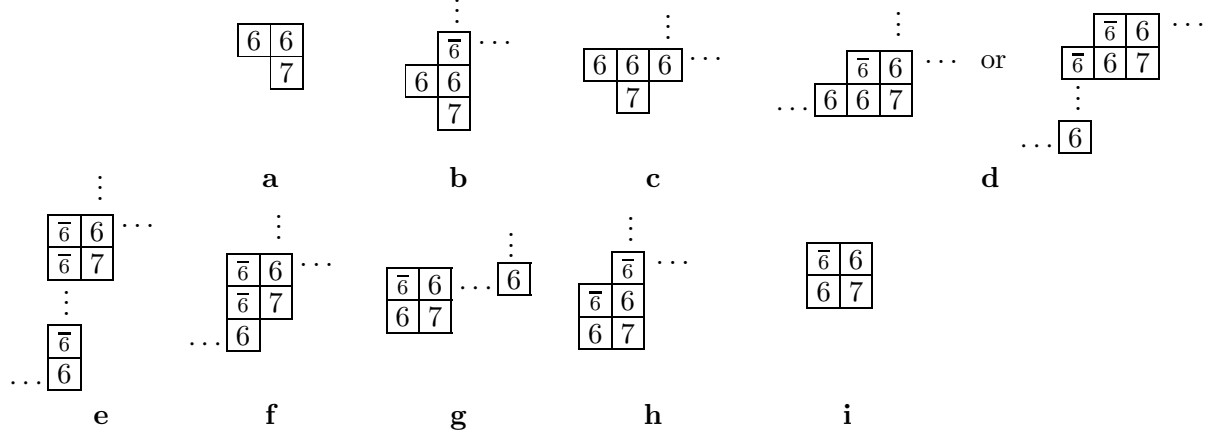
\vspace{-7mm}
$$
\begin{array}{cccccc}
&
\begin{array}{l}
\young(66,:7)
\end{array}
&
\begin{array}{rl}
\vdots\,\,\vspace{-1.7mm}\\
&\!\!\!\!\ldots\vspace{-.3cm}\\
\young(:\tsi,66,:7)
\end{array}
&
\begin{array}{rl}
\vdots\,\,\vspace{-1.7mm}\\
&\!\!\!\!\ldots\vspace{-.3cm}\\
\young(666,:7)
\end{array}
&

\begin{array}{rl}
\vdots\,\,\vspace{-1mm}\\
&\!\!\!\ldots\vspace{-.3cm}\\
\cdots\young(:\tsi6,667)
\end{array}

\textrm{ or}

\begin{array}{rl}
\vdots\,\,\vspace{-1mm}\\
&\!\!\!\ldots\vspace{-.3cm}\\
\young(:\tsi6,\tsi67)\vspace{-2.5mm}\\
\vdots\ \ \ \ \ \ \ \ \\
\cdots\young(6)\ \ \ \ \ \ \,\, 
\end{array}

\\
&\textrm{\bf a}&\textrm{\bf b}&\textrm{\bf c}&
\textrm{\bf d}\vspace{-5mm}\\

\begin{array}{ll}
\ \ \ \ \ \vdots\vspace{-1mm}\\
&\!\!\!\!\ldots\vspace{-.3cm}\\
\young(\tsi6,\tsi7)\vspace{-2mm}\\
\,\,\,\vdots\\
 \!\!\!\!\!\!\!\!\cdots\young(\tsi,6)
\end{array}
&
\begin{array}{ll}
\ \ \ \ \ \ \ \vdots\vspace{-1mm}\\
&\!\!\!\!\ldots\vspace{-.3cm}\\
\!\!\!\!\cdots\young(\tsi6,\tsi7,6)
\end{array}
&
\begin{array}{rl}
&\ \ \ \vdots\vspace{-.7mm}\\
&\!\!\!\!\cdots\young(6)\vspace{-5mm}\\
 \young(\tsi6,67)
\end{array}
&
\hspace{-.5cm}
\begin{array}{l}
\ \ \ \ \ \vdots\vspace{-1.5mm}\\
\ \ \ \ \ \ \ \ \ldots\vspace{-3mm}\\
\young(:\tsi,\tsi6,67)
\end{array}
&
\hspace{-4cm}
\begin{array}{c}
\young(\tsi6,67)
\end{array}\\
\textrm{\bf e}&
\textrm{\bf f}&
\textrm{\bf g}&\textrm{\bf h}&\hspace{-3.8cm}\textrm{\bf i}
\end{array}
\vspace{-1mm}$$
Note that in the cases d,e,f,g and h given above we can modify the
filling to get another amenable one as shown in the following figures.
Therefore if $P_m\cup P_{m-1}$ is one of those cases
then $\DD$
cannot be strange.
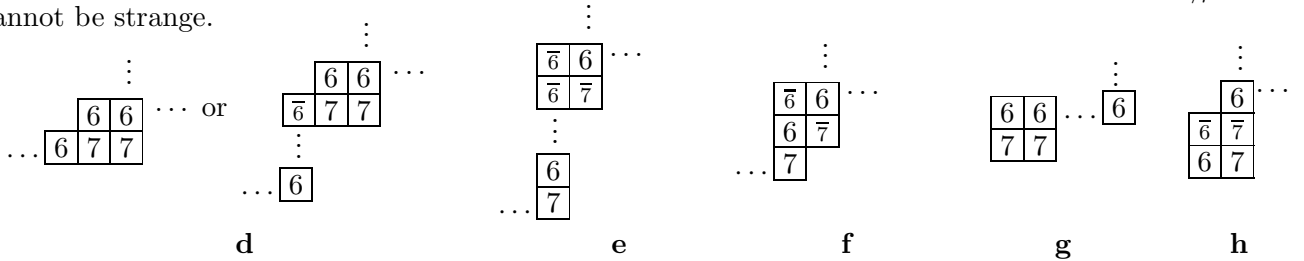
\vspace{-.8cm}
$$
\begin{array}{ccccc}


\begin{array}{rl}
\vdots\,\,\vspace{-1mm}\\
&\!\!\!\ldots\vspace{-.3cm}\\
\cdots\young(:66,677)
\end{array}

\textrm{or}

\begin{array}{rl}
\vdots\,\,\vspace{-1mm}\\
&\!\!\!\ldots\vspace{-.3cm}\\
\young(:66,\tsi77)\vspace{-2.5mm}\\
\vdots\ \ \ \ \ \ \ \ \\
\cdots\young(6)\ \ \ \ \ \ \,\, 
\end{array}

\qquad &


\begin{array}{ll}
\ \ \ \ \ \vdots\vspace{-1mm}\\
&\!\!\!\!\ldots\vspace{-.3cm}\\
\young(\tsi6,\tsi\tse)\vspace{-2mm}\\
\,\,\,\vdots\\
 \!\!\!\!\!\!\!\!\cdots\young(6,7)
\end{array}
\qquad
&


\begin{array}{ll}
\ \ \ \ \ \ \ \vdots\vspace{-1mm}\\
&\!\!\!\!\ldots\vspace{-.3cm}\\
\!\!\!\!\cdots\young(\tsi6,6\tse,7)
\end{array}

\qquad &


\begin{array}{rl}
&\ \ \ \vdots\vspace{-.7mm}\\
&\!\!\!\!\cdots\young(6)\vspace{-5mm}\\
 \young(66,77)
\end{array}
&
\begin{array}{l}
\ \ \ \ \ \vdots\vspace{-1.5mm}\\
\ \ \ \ \ \ \ \ldots\vspace{-3mm}\\
\young(:6,\tsi\tse,67)
\end{array}\\
\textrm{\bf d}&\textrm{\bf e}&\textrm{\bf f}&\textrm{\bf g}&
\textrm{\bf h}
\end{array}
$$


The only remaining possibility is Case i. 
If $m=2$ then $\DD$ is equal to the underlying diagram of 
figure i, and 
has a unique amenable filling. Assume $m>2$. 
Next we prove Proposition \ref{casei}.

\begin{definition}
A path of boxes is called a $(p,q)$-hook if it has $p$ vertical 
boxes in the first column and $q$ horizontal boxes 
in the first row.
\end{definition}

\begin{example}\rm
The following figure illustrates a $(2,3)$-hook.
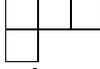
\vspace{-1mm}
$$
\young(\ \ \ ,\ )
\vspace{-4mm}
$$
\end{example}

\begin{proposition}\label{casei}
Suppose $\DD$ is a strange diagram and
$P_{m-1}\cup P_m$ is as in Case i. 
Then for any $j$, $P_j$ 
is a $(p_j,m-j+1)$-hook
for some $p_j\leq m-j+1$. Moreover, if 
for some $j>1$ we have $p_j<m-j+1$, then $p_j=p_{j-1}=\cdots=p_1$.
\vspace{-2mm}
\end{proposition}

\begin{proof}
We use backwards induction. Assume the statement holds for
$$P_{m-r+1},...,P_m$$ where $r>1$. 
We prove it for $P_{m-r}$. 
Suppose that $P_{m-r+1}$ is an $(l,r)$-hook 
for some $1\leq l\leq r$. 
As before, assume 
$\BB_{x,y}\in P_m$.
Lemma \ref{discon}, Lemma \ref{paths} 
and
Lemma \ref{propD} imply that
\item[$\diamond$]
For any $j$ such that $1\leq j\leq l-1$, 
we have
$\BB_{x-r,y+r-j}\in P_{m-r}$

\noindent and
\item[$\diamond$] For any $j$ such that 
$0\leq j\leq r$, we have $\BB_{x-j,y+r}
\in P_{m-r}$.

Let $Q$ denote
the path 
formed by the boxes
in the two statements given above.
Obviously $Q\subseteq P_{m-r}$.
We have the following two cases:

\begin{itemize}
\item[\bf Case A.] $l<r$. 
Then $\BB_{x-r,y+r+1-l}$ should be the last box of
$P_{m-r}$.
 We will show that $P_{m-r}=Q$. 
Suppose, on the contrary, that $P_{m-r}\neq Q$.
Then either $\BB_{x,y+r+1}\in P_{m-r}$
or $\BB_{x+1,y+r}\in P_{m-r}$. However, 
$\DD$ cannot be strange in 
either of the cases, as shown below.

\begin{itemize}

\item[\bf Case A1.] $\BB_{x+1,y+r}\in P_{m-r}$. Let 
$u=x-r,v=y+r+1-l$. Then the following 
procedure provides a new filling. 
\begin{itemize}
\item[\tt Step 1.] Set $t=1$ and $\BB_{u,v+1}\leftarrow m-r$.
\item[\tt Step 2.] If $\BB_{u+t-1,v-t+1}\notin\DD$ then 
stop.
\item[\tt Step 3.] Find smallest $y'\geq v-t$ such that
$\BB_{u+t,y'}\in\DD$. If there is no such $y'$, then stop.
\item[\tt Step 4.] If $\BB_{u+t,y'}\in P_s$ then 
do $\BB_{u+t-1,v-t+1}\leftarrow s$ and 
$\BB_{u+t,v-t+1}\leftarrow s$.
\item[\tt Step 5.] Increase $t$ by one. Go back to Step 2.
\end{itemize}
The idea behind
the procedure can be more
concretely described as follows. We eliminate the last box of
$P_{m-r}$, add this box to $\DD^{(m-r+1)}$, and find new paths
$P_{m-r+1},...,P_m$ according to the algorithm of section
\ref{algorithm}.

The following figure demonstrates 
the above procedure when $m=8$
and $r=3$.
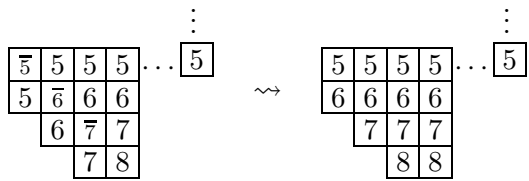
\vspace{-4mm}
$$
\begin{array}{ccc}
\begin{array}{rl}
&\ \ \,\vdots\\
&\!\!\!\!\!\cdots\young(5)\vspace{-5mm}\\
\young(\tfi555,5\tsi66,:6\tse7,::78)
\end{array}
&
\leadsto
&
\begin{array}{rl}
&\ \ \,\vdots\\
&\!\!\!\!\!\cdots\young(5)\vspace{-5mm}\\
\young(5555,6666,:777,::88)
\end{array}
\end{array}
\vspace{-2mm}
$$

\item[\bf Case A2.]
$\BB_{x,y+r+1}\in P_{m-r}$. Then the following procedure
gives the new amenable
filling. 
\begin{itemize}
\item[\tt Step 1.] Set $t=r+1$.
\item[\tt Step 2.] If $t>1$ then 
replace the entry of $\BB_{x,y+t}$
by that of $\BB_{x,y+t-1}$. If $t=1$ then replace the 
entry of $\BB_{x,y+t}$ by $\ov{m}$. If $t<1$ then stop.
\item[\tt Step 3.] Decrease $t$ by one. Go back to Step 2.\\
\vspace{-3mm}
\end{itemize}
The amenability of the diagram obtained by the procedure 
follows from statements similar to those that appear 
in the proof of Lemma \ref{amenamen}.

For $m=7$ and $r=2$ the following figure demonstrates 
the effect of the procedure.
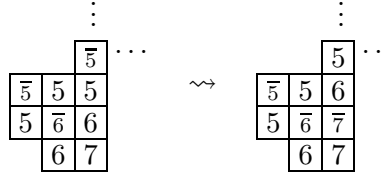
\vspace{-3.5mm}
$$
\begin{array}{ccc}
\begin{array}{rl}
\vdots\ \vspace{-1mm}\\
&\!\!\!\!\ldots\vspace{-3mm}\\
\young(::\tfi,\tfi55,5\tsi6,:67)\\
 \end{array}
&
\leadsto
&
\begin{array}{rl}
\vdots\ \vspace{-1mm}\\
&\!\!\!\!\ldots\vspace{-3mm}\\
\young(::5,\tfi56,5\tsi\tse,:67)\\
\end{array}\\
\end{array}
\vspace{-2mm}
$$

\end{itemize}

\item[\bf Case B.] 
$l=r$. 
In this case if $\BB_{x-r,y}\notin P_{m-r}$
then we can argue as in the case $l<r$ to show that
$P_{m-r}=Q$. The following figure illustrates what
happens when $m=7$, $r=2$ and $\BB_{x+1,y+r}\in P_{m-r}$.
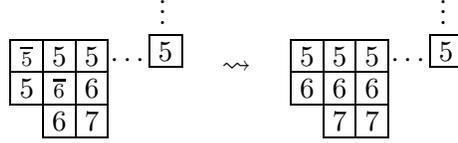
\vspace{-3.5mm}
$$
\begin{array}{ccc}
\begin{array}{rl}
&\ \ \,\vdots\\
&\!\!\!\!\!\cdots\young(5)\vspace{-5mm}\\
\young(\tfi55,5\tsi6,:67)
\end{array}
&
\leadsto
&
\begin{array}{rl}
&\ \ \,\vdots\\
&\!\!\!\!\!\cdots\young(5)\vspace{-5mm}\\
\young(555,666,:77)
\end{array}
\end{array}
\vspace{-2mm}
$$

Therefore we may assume 
$\BB_{x-r,y}\in P_{m-r}$. Next we 
show that $$P_{m-r}=Q\cup\{\BB_{x-r,y}\}.$$
Suppose on the contrary that  
$
P_{m-r}\neq Q\cup\{
\BB_{x-r,y}\}$.
Then  
one of the following cases occurs. In all of them we show that $\DD$ 
is not strange.

\begin{itemize}

\item[\bf Case B1.] $\BB_{x-r-1,y}\in P_{m-r}$.
We can get a new amenable filling by 
the operations $\BB_{x-r,y}\leftarrow m-r+1$ 
and 
$\BB_{x-r,y+1}\leftarrow m-r$.

\item[\bf Case B2.] $\BB_{x-r,y-1}\in P_{m-r}$.
In this case we can obtain a new amenable filling as follows.
Find the smallest $y'$ such that $\BB_{x-r,y'}\in P_{m-r}$. Then 
do $\BB_{x-r,y'+1}\leftarrow m-r$, $\BB_{x-r,y'}\leftarrow m$
and $\BB_{x,y}\leftarrow\ov{m}$.

\item[\bf Case B3.]$\BB_{x-r-1,y}\notin P_{m-r}$, 
$\BB_{x-r,y-1}\notin P_{m-r}$ 
but
$\BB_{x+1,y+r}\in P_{m-r}$.
The new amenable filling is obtained by 
$$\BB_{x-r,y}\leftarrow m-r+1\textrm{ and }
\BB_{x-r,y+1}\leftarrow m-r.$$
The following figure illustrates this procedure when $m=7,r=4$.
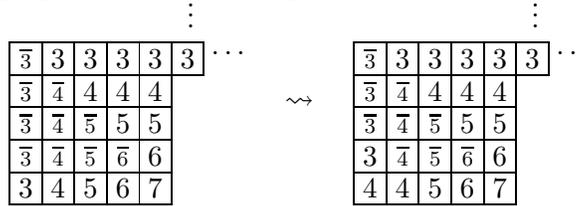
\vspace{-3.5mm}
$$
\begin{array}{ccc}
\begin{array}{rl}
\vdots\ \vspace{-1mm}\\
&\!\!\!\!\ldots\vspace{-3mm}\\
\young(\tth33333,\tth\tfo444,\tth\tfo\tfi55,\tth\tfo\tfi\tsi6,34567)\\
\end{array}
&
\leadsto
&
\begin{array}{rl}
\vdots\ \vspace{-1mm}\\
&\!\!\!\!\ldots\vspace{-3mm}\\
\young(\tth33333,\tth\tfo444,\tth\tfo\tfi55,3\tfo\tfi\tsi6,44567)\\
\end{array}
\end{array}
\vspace{-2mm}
$$
\item[\bf Case B4.] 
$\BB_{x-r-1,y}\notin P_{m-r}$,
$\BB_{x-r,y-1}\notin P_{m-r}$ but
$\BB_{x,y+r+1}\in P_{m-r}$.
In this case exactly the same 
procedure that was
given in the analysis of Case
A2 can provide a new filling. We illustrate it
below with $m=7$ and $r=4$.
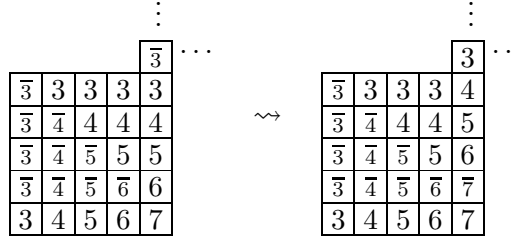
\vspace{-3mm}
$$
\begin{array}{ccc}
\begin{array}{rl}
\vdots\ \vspace{-1mm}\\
&\!\!\!\!\ldots\vspace{-3mm}\\
\young(::::\tth,\tth3333,\tth\tfo444,\tth\tfo\tfi55,\tth\tfo\tfi\tsi6,34567)

\end{array}
&
\leadsto
&
\begin{array}{rl}
\vdots\ \vspace{-1mm}\\
&\!\!\!\!\ldots\vspace{-3mm}\\
\young(::::3,\tth3334,\tth\tfo445,\tth\tfo\tfi56,\tth\tfo\tfi\tsi\tse,34567)
\end{array}
\end{array}
\vspace{-7mm}
$$
\end{itemize}
\end{itemize}       
\end{proof}
\noindent{\bf Remark.} The proof of
Proposition \ref{casei}
actually proves more. It proves the following corollary.

\begin{corollary}\label{caseicor}
Let $\DD$ be a strange diagram and
$P_m\cup P_{m-1}$ be as in Case i. Then 
$\lambda/\mu=\overline{\lambda}/\overline{\mu}$ such that 
\begin{enumerate}
\item[$\bullet$] $\overline{\lambda}=\{p+q+r>p+q+r-1>p+q+r-2>\cdots>p\}$ and 
$\overline{\mu}=\{q>q-1>\cdots>1\}$ where $p,q,r$ are integers such 
that $p,q\geq 1$, $r\geq 0$. 
\item[$\bullet$] $\overline{\lambda}=\{p+q>p+q-1>\cdots>p+q-r\}$ and 
$\overline{\mu}=\{q>q-1>\cdots>q-r\}$
where $p,q,r$ are integers such that $p>0$ and $q>r\geq 0$.
\end{enumerate}

\end{corollary}

\noindent{\bf Remark.}
The two cases of Corollary \ref{caseicor}
may overlap. Moreover, not all of them are such that 
$P_m\cup P_{m-1}$ is as in Case i. However, we have written 
it in the given form above so that we can refer to it at other points of
this manuscript. For example see Corollary \ref{IIcor}.

Next we analyze cases a,b and c which were introduced at the beginning
of section \ref{CASE_I}.
\vspace{-2mm}
\begin{proposition}\label{caseb}
If $\DD$ is a strange diagram and
$P_{m}\cup P_{m-1}$ is as in Case b then
we have 
$$\lambda=\{m>m-1>\cdots>1\}\textrm{ and }\mu=\{\mu_1>\cdots>\mu_l\}
\textrm{ where }l<m-1.$$ 
\end{proposition}
\begin{proof}

It suffices to show that there does not exist any $z$ such that
$\BB_{x+1,z}\in\DD$. Assume the contrary, and let $z$ be the smallest
integer such that $\BB_{x+1,z}\in\DD$. 
Let $r$ be such that $\BB_{x+1,z}\in P_r$.
From Lemma \ref{propD}, Lemma \ref{propT}
and Lemma \ref{discon}
it follows that $\BB_{x,z}\in P_r$.
Moreover, by Lemma \ref{paths} it follows that for any $s$ such that
$r\leq s<m$, 
$$
\BB_{x-m+s+1,y+m-s+1},\BB_{x-m+s+1,y+m-s},\BB_{x-m+s,y+m-s}\in P_s.
$$ 
Taking $s=r$ implies that  
$z>y+m-r$ and  
the path $P_r$
makes a clockwise turn at some point.

Let $u,v$ be such that $\BB_{u,v},\BB_{u-1,v},\BB_{u,v+1}\in P_r$
and $v<z$. For example we can take $u=x-m+r+1,v=y+m-r$.
It follows from Lemma \ref{paths} that for any $t>0$ such that $t+r<m$, 
$$
\BB_{u+t,v-t},\BB_{u+t-1,v-t},\BB_{u+t,v-t+1}\in P_{r+t}.
$$
Now the following procedure provides a new amenable filling.

\begin{itemize}
\item[\tt Step 1.] Set $t=0$.
\item[\tt Step 2.] If $u+t<x$ then 
do 
$\BB_{u+t,v-t}\leftarrow \ov{r+t+1}$ and 
$\BB_{u+t+1,v-t}\leftarrow r+t+1$.
If $u+t=x$ and the entry of $\BB_{u+t,v-t-1}$  
is either $w$ or $\ov{w}$ where $w\in\mathcal A$, 
then do
$\BB_{u+t,v-t}\leftarrow\ov{w}$. If $u+t>x$ then stop.
\item[\tt Step 3.] Increase $t$ by one. Go back to Step 2.
\end{itemize}
The idea of the procedure is more concretely explained as follows.
We eliminate one of the boxes of $P_r$, add it to $\DD^{(r+1)}$,
and then find new paths $P_{r+1},...,P_m$ according to
the algorithm of section \ref{algorithm}.

The following figure demonstrates how the procedure can be applied
to a case with $m=7$ and $r=4$. In fact it can be seen that there
are at least two ways to get new amenable fillings.
\vspace{-3mm}
$$
\begin{array}{ccccc}
\begin{array}{rl}
\vdots\ \vspace{-1mm}\\
&\!\!\!\!\ldots\vspace{-3mm}\\
\young(::\tfo44,:\tfo4\tfi,:\tfo\tfi5,44\tfi\tsi,:55\tsi,::66,:::7)
\end{array}
&
\leadsto
&
\begin{array}{rl}
\vdots\ \vspace{-1mm}\\
&\!\!\!\!\ldots\vspace{-3mm}\\
\young(::\tfo44,:\tfo4\tfi,:\tfo\tfi5,4\tfi5\tsi,:5\tsi6,::6\tse,:::7)
\end{array}
&
\textrm{or}
&
\begin{array}{rl}
\vdots\ \vspace{-1mm}\\
&\!\!\!\!\ldots\vspace{-3mm}\\

\young(::\tfo44,:\tfo\tfi5,:\tfo\tfi\tsi,44\tfi\tsi,:55\tsi,::66,:::7)

\end{array}
\end{array}
\vspace{-7mm}
$$
\end{proof}
\begin{proposition}\label{casea}
Let $\DD$ be a strange diagram. 
Assume $P_m\cup P_{m-1}$ is as in Case a or Case c. 
Then
$\DD$ is given by one of the following cases.

\begin{enumerate}
\item[$\bullet$]
$\DD$ is given by partitions of Proposition \ref{caseb}.
\item[$\bullet$] $\DD=D_{\overline{\lambda}}$ 
for an arbitrary $\overline{\lambda}$. 
\end{enumerate}
\end{proposition}

\begin{proof}
Let, as before, $\BB_{x,y}\in P_m$.
Suppose $\DD$ is not given by partitions of Proposition
\ref{caseb}. According to the proof of 
Proposition \ref{caseb},
this means that there exists a $t$ such that 
$\BB_{x+1,t}\in\DD$. We can assume $t$ is the smallest integer 
with this property. Let $r$ be such that $\BB_{x+1,t}\in P_r$.
Let $r'$ be the largest integer for which $P_{r'}$ does not 
lie within a single row. Note that if there is no such $r'$ then 
it follows that $\DD=D_{\lambda}$.

Set $r_1=\min\{r,r'\}$. An argument 
similar to that of Proposition \ref{caseb} shows that $P_{r_1}$
should make a clockwise turn; i.e. there exist $p,q$
such that $\BB_{p,q},\BB_{p-1,q},\BB_{p,q+1}\in P_{r_1}$.
Now one of the following possibilities happens.
\begin{itemize}
\item[\bf Case C1.] $r>r_1$. Then it follows that $r>r'$. In 
this case we use the following procedure.
\begin{enumerate}
\item[\tt Step 1.] Set $t=0$. 
\item[\tt Step 2.] If $\BB_{p,q+t}\in\DD$ then do 
$\BB_{p,q+t+1}\leftarrow\BB_{p,q+t}$. Otherwise,
do $\BB_{p,q+t+1}\leftarrow r_1-t$ and stop.
\item[\tt Step 3.] Decrease $t$ by one. Go back to Step 2.\\
\end{enumerate}
The idea of the procedure is similar to that of the
procedure given in Case B4 of Proposition \ref{casei}.
\item[\bf Case C2.] $r=r_1$. Then it follows that
$r\leq r'$. In this case we use the following procedure.
\begin{itemize}
\item[\tt Step 1.] Set $t=1$ and do $\BB_{p,q+1}\leftarrow r_1$.
\item[\tt Step 2.] If $\BB_{p+t,q-t}\in \DD$, then 
do $\BB_{p+t-1,q-t+1}\leftarrow \ov{r_1+t}$ and 
$\BB_{p+t,q-t+1}\leftarrow r_1+t$.
\item[\tt Step 3.] If $\BB_{p+t,q-t}\notin\DD$ 
and $\BB_{p+t-1,q-t}\in \DD$ then
do $\BB_{p+t-1,q-t+1}\leftarrow \ov{r_1+t}$ and 
stop.
\item[\tt Step 4.] If $\BB_{p+t,q-t},\BB_{p+t-1,q-t}\notin\DD$
then do $\BB_{p+t-1,q-t+1}\leftarrow r_1+t$ and stop.
\item[\tt Step 5.] Increase $t$ by one. Go back to Step 2.
\end{itemize}
\vspace{1mm}
The idea behind the procedure is to eliminate one of the boxes of
$P_{r_1}$, modify the filling of $P_{r_1}$ slightly, attach the 
box removed
from $P_{r_1}$ to $\DD^{(r_1+1)}$ and 
find new paths $P_{r_1+1},...,P_m$ according to the 
algorithm of section \ref{algorithm}.
\end{itemize}

\end{proof}

\subsection{Case II: $P_m$ lies within a single
row and has at least two boxes}
Let us assume that the boxes that belong to $P_m$ are
$\BB_{x-l,y},...,\BB_{x,y}$ where $l>0$. From Lemma 
\ref{propT} it follows that for any $j\in\{0,...,l+1\}$,
$\BB_{x-j,y+1}\in P_{m-1}$.
We have the following two cases.

\vspace{1mm}
{\bf Case D1.} $\BB_{x-l-1,y}\in\DD$. In this case the boxes  
$\BB_{x-l-1,y+1},...,\BB_{x,y+1}$ and
$\BB_{x-l-1,y}$ form a path in $P_{m-1}$. 
An argument similar to the one given in
Case B of Proposition \ref{casei}
proves that there are no other boxes in $P_{m-1}$.
Next we state the
following proposition which 
is similar in statement and proof to Proposition 
\ref{casei}. 
\begin{proposition}
\label{case1row}
Suppose that $\DD$ is a strange diagram such that
$P_m\cup P_{m-1}$ is as in Case D1. Then 
for any 
$j\in\{1,...,m\}$, $P_j$ is a $(p_j,l+m-j+1)$-hook
for some $p_j\leq m-j+1$. 
Moreover, if there exists a $j>1$ such that $p_j<m-j+1$, then
$p_j=p_{j-1}=\cdots=p_1$. 
\end{proposition}
We omit the proof of Proposition \ref{case1row} because it is very 
similar to the proof of Proposition \ref{casei}.

\vspace{1mm}
{\bf Case D2.}
$\BB_{x-l-1,y}\notin\DD$.
The analysis of this case is similar to that of Case a in section 
\ref{CASE_I}. (See Proposition \ref{casea}.) 
In fact we can prove the following proposition.
\begin{proposition}
Suppose that $\DD$ is a strange diagram which is as 
in Case II above,
and moreover $\BB_{x-l-1,y}\notin\DD$. Then 
$\DD=D_{\overline{\lambda}}$ for an arbitrary $\overline{\lambda}$.
\end{proposition}
\begin{proof}
It suffices to show that for any $j$, $P_j$ lies within a
single row. The proof is by contradiction. 
Assume the contrary, and let $j$ be the largest integer such that
$P_j$ does not lie within one row. From Lemma \ref{propT}
it follows that if $p'$ is the largest integer for which 
there exists a box $\BB_{p',y+m-j-1}\in P_{j+1}$ then 
there exists a $p\geq p'$ such that 
$\BB_{p,q},\BB_{p-1,q},\BB_{p,q+1}\in P_j$, where
$q=y+m-j$. 
Now set $r_1=j$ and 
apply the procedure given in Case C1 of the proof of Proposition 
\ref{casea}
to get a new amenable filling.
\vspace{-3mm}
\end{proof}
\begin{corollary}\label{IIcor}
Let $\DD$ be as in Case II. Suppose $\DD$ is a strange diagram.
Then either $\DD=D_{\overline{\lambda}}$ for an arbitrary $\overline{\lambda}$,
or $\DD$ is given as in
the statement of Corollary \ref{caseicor}.
\end{corollary}

\subsection{Case III: $P_m$ lies within a single column
and has at least two boxes}
The analysis in this case is pretty similar to the 
previous cases. 
Suppose $P_m$ consists of the boxes $\BB_{x,y+l},...,\BB_{x,y}$
where $l>0$.
Then for any $j\in\{1,...,l+1\}$, $\BB_{x-1,y+j}\in P_{m-1}$.
Moreover, $\BB_{x,y+l+1}\in P_{m-1}$. 
Next we show that none of the boxes
$\BB_{x+1,y+l+1},\BB_{x-1,y-1},\BB_{x-2,y}$ can 
belong to $P_{m-1}$. In fact in each of 
the following cases we show that $\DD$ cannot be a strange diagram.

\begin{itemize}
\item[$\bullet$] $\BB_{x-1,y-1}\in P_{m-1}$. Then we get a new 
filling as follows. Find the smallest $y'$ such that
$\BB_{x-1,y'}\in P_{m-1}$ and do 
\vspace{-2mm}
$$\BB_{x-1,y'+1}\leftarrow m-1,
\BB_{x-1,y'}\leftarrow m\textrm{ and }\BB_{x,y}\leftarrow \ov{m}.
\vspace{-1.5mm}$$
\item[$\bullet$] $\BB_{x-2,y}\in P_{m-1}$. Then we obtain a new
filling by 
\vspace{-2mm}
$$\BB_{x-1,y}\leftarrow m\textrm{ and }
\BB_{x-1,y+1}\leftarrow 
m-1.
\vspace{-1.5mm}$$     
\item[$\bullet$] $\BB_{x-1,y-1},\BB_{x-2,y}\notin P_{m-1}$ 
but $\BB_{x+1,y+l+1}\in P_{m-1}$. Then a new filling can be obtained
as in
Case A1 or Case B3 of Proposition \ref{casei}. We leave the details to
the reader.
\end{itemize}
\begin{proposition}\label{verticalI}
Let $\DD$ be a strange diagram. Suppose 
$\BB_{x-1,y}\notin P_{m-1}$.
Then 
$$\lambda=\{m>m-1>\cdots>1\} \textrm{ and } 
\mu=\{\mu_1>\cdots>\mu_l\}$$
where $l<m-1$. 
\end{proposition}
\begin{proof}
It suffices to show that there does not exist any $y'$ such that 
$\BB_{x+1,y'}\in\DD$. Suppose there exists such a $y'$, and without loss 
of generality assume that $y'$ is the smallest integer with this property.
Suppose $\BB_{x,y'}\in P_r$. Then by Lemmas \ref{paths} and \ref{propT} 
it follows that 
\vspace{-2mm}$$
\BB_{x-m+r,y+m-r},\BB_{x-m+r,y+m-r+1}\in P_r.
\vspace{-1mm}$$
Now the following procedure provides a new amenable filling.
\begin{enumerate}
\item[\tt Step 1.]
Set $j=r$ and do $\BB_{x-m+r,y+m-r+1}\leftarrow r$.
\item[\tt Step 2.]
If $j<m$ then do $\BB_{x-m+j,y+m-j}\leftarrow j+1$ and 
$\BB_{x-m+j+1,y+m-j}\leftarrow j+1$.
\item[\tt Step 3.] If $j=m$ then do $\BB_{x-m+j,y+m-j}\leftarrow m+1$
and stop.
\item[\tt Step 4.] Increase $j$ by one. Go back to Step 2.
\end{enumerate}
\vspace{-7mm}
\end{proof}
The proof of the following proposition is very similar to those of
Propositions \ref{case1row} and \ref{casei}, 
and therefore we omit its proof.
\begin{proposition}\label{verticalII}
Let $\DD$ be a strange diagram given as in
Case III and suppose 
$\BB_{x-1,y}\in P_{m-1}$. Then $\lambda,\mu$ are
given as in Corollary \ref{caseicor}.
\end{proposition}

\begin{corollary}
Let $\DD$ be a strange diagram given as in Case III. Then 
$\DD$ is given either by the partitions of
Corollary 
\ref{caseicor}
or by the partitions of Proposition
\ref{verticalI}.
\end{corollary}

\subsection{Proof of uniqueness} In the previous sections we 
showed that the only diagrams which could possibly be strange 
are those
listed in Theorem \ref{main}.
In this section we prove that all of those diagrams are indeed 
strange. First we need a simple property of 
any arbitrary GSYT.

\begin{lemma}\label{strict}
Let $T$ be an amenable GSYT and let $w(T)$, the row word of $T$, have 
length $n$. Then for any $k>1$ such that $m_{k-1}(n)>0$ 
we have $m_k(n)<m_{k-1}(n)$.
\end{lemma}
\begin{proof}
From Definition \ref{amenable}
it follows that for any $k>1$, $m_k(n)\leq m_{k-1}(n)$.
Suppose $m_{k-1}(n)=m_{k}(n)$. Then it follows 
from Definition \ref{amenable}
that 
if $j_k$ (respectively $j_{k-1}$)
is the smallest integer such that $w_{j_k}=k$ 
(respectively $w_{j_{k-1}}=k-1$) then $j_k<j_{k-1}$.
The fourth part of Definition \ref{amenable}
implies that $w_j\neq\ov{k-1}$ for any $j<j_{k-1}$.
Since $m_{k-1}(j)\geq m_k(j)$ for any $j\in\{0,...,2n-1\}$,
it follows that
$m_k(n+j_{k-1}-1)=m_{k-1}(n+j_{k-1}-1)=m_{k-1}(n)$.
But then
$w_{j_{k-1}}=k-1$, which contradicts the second part of
Definition \ref{amenable}.
\vspace{-3mm}
\end{proof}
It is obvious that if 
$\DD=D_\lambda$ then $\DD$ is strange. Next we prove that 
if $\DD$ is given by partitions of Corollary
\ref{caseicor} then $\DD$ is a strange diagram. We will 
present the argument only for the case  
$\lambda=\{p+q+r>p+q+r-1>p+q+r-2>\cdots>p\}$ and 
$\mu=\{q>q-1>\cdots>1\}$ where $p,q,r$ are integers such 
that $p,q\geq 1$, $r\geq 0$. The other case can be treated
in a very similar fashion.
\begin{proposition}\label{semirec0}
Let $\DD$ be given by 
$\lambda=\{p+q+r>p+q+r-1>p+q+r-2>\cdots>p\}$ and 
$\mu=\{q>q-1>\cdots>1\}$ where $p,q,r$ are integers such 
that $p,q\geq 1$, $r\geq 0$. Then $\DD$ is a strange diagram.
\end{proposition}
We give the proof of Proposition \ref{semirec0} 
through Lemma \ref{semirec1}, Lemma \ref{semirec2} and Lemma
\ref{semirec3}.
Let us assume that $\BB_{x,y}$
is the top right box of the diagram $\DD$. 
Let $n$ denote the number of boxes of $\DD$.
We consider an arbitrary 
amenable
filling of $\DD$ and show that it has to be the one obtained by
the algorithm of section \ref{algorithm}.
\vspace{-1.5mm}
\begin{lemma}\label{semirec1}
For any integers 
$k,x'$ such that $k\geq 0$ and 
$x-p-r+2+k\leq x'\leq x$, if 
$\BB_{x',y-k}\in\DD$ then it
should be filled by a $k+1$.
\end{lemma}

\begin{proof}
We prove this lemma by induction on $k$. For $k=0$,
the above statement 
follows from the fact that we should have 
$\BB_{x,y}\in\{\ov{1},1\}$ and moreover
no two $\ov{1}$'s can lie within 
the same row. Next assume $k\geq 1$. Let $x'$ be chosen as 
above. Then by induction hypothesis, the boxes
$\BB_{x',y-k+1},\BB_{x'-1,y-k+1}$ are filled by $k$.
Moreover, none of the rows which lie above the box
$\BB_{x,y-k}$
can contain an element of $\mathcal{A}$ which is
strictly larger than $k$, 
because their rightmost 
boxes are filled by elements less than or equal to $k$.
Consequently, $\BB_{x,y-k}$ should be filled by
an element of $\{\ov{k+1},k+1\}$. (Note that $\BB_{x,y-k+1}$
contains a $k$.) This in turn limits the  
entries of $\BB_{x',y-k}$ and $\BB_{x'-1,y-k}$ 
to $\{\ov{k+1},k+1\}$. However, if
$\BB_{x',y-k}$ is filled by a $\ov{k+1}$ 
then $\BB_{x'-1,y-k}$ should also be filled by a $\ov{k+1}$,
which is a contradiction by
Definition \ref{GSYT}.
Therefore $\BB_{x',y-k}$ has to be filled by a $k+1$.
\end{proof}

Lemma \ref{semirec1} implies that for any $k\geq 1$,
if there is at least one box in $\DD$ filled by a $k$, 
then there are at
least $p+r-k$ boxes in $\DD$ which are filled by a $k$. Note that 
since $\DD$ has $p+r$ columns, there can be at most $p+r$ 
boxes which are filled by a $1$. Lemma \ref{strict} implies 
that there can be at most $p+r-k+1$ boxes which are 
filled by a $k$. We will see below that actually the latter
possibility takes place. More accurately, 
we show that if there is at least 
one $k$ in the filling of $\DD$, then there are exactly
$p+r-k+1$ of them.
\vspace{-1mm}
\begin{lemma}\label{semirec2}
The largest integer that
appears in the filling of $\DD$ is $\min\{p+r,q+r+1\}$.
For any $k\in\{1,...,\min\{p+r,q+r+1\}\}$
there are exactly $p+r-k+1$ boxes in $\DD$
which contain a $k$.
\end{lemma}

\begin{proof}
If $p+r>q+r+1$ then by Lemma \ref{semirec1} 
the box $\BB_{x,y-q-r}$
is filled by a $q+r+1$. But $\BB_{x,y-q-r}$ is the rightmost
box of the lowest row of $\DD$, and 
since the rows and columns of a filling are 
weakly increasing, all of the boxes
have to be filled by elements of $\mathcal{A}$ which are less than 
or equal to $q+r+1$. Now 
Lemma \ref{semirec1} implies that $\BB_{x+q-p+2,y-q-r}$ contains a 
$q+r+1$ and $\BB_{x+q-p+1,y-q-r+1}$ contains a $q+r$. Because
of the shape of $\DD$ 
we can see that $\BB_{x+q-p+1,y-q-r}\in\DD$.
Since $\BB_{x+q-p+1,y-q-r+1}$ contains a $q+r$, 
the entry of $\BB_{x+q-p+1,y-q-r}$ 
should belong to $\{\ov{q+r+1},q+r+1\}$. But 
the entry of $\BB_{x+q-p+1,y-q-r}$ cannot be
$\ov{q+r+1}$ because this implies that the 
leftmost occurence
of an element of
$\{\ov{q+r+1},q+r+1\}$
in the 
row word of the filling is marked.
Consequently, 
$\BB_{x+q-p+1,y-q-r}$ is filled by a 
$q+r+1$ 
which implies that
$m_{q+r+1}(n)\geq p-q$. Now we have
$$
p+r\geq m_1(n)>\cdots>m_{q+r+1}(n)\geq p-q
$$
which implies that $m_k(n)=p+r-k+1$ for any $k\in\{1,...,q+r+1\}$.

If $p+r\leq q+r+1$, then Lemma 
\ref{semirec1} implies that $\BB_{x,y-p-r+2}$ 
has to be filled by $p+r-1$. Now $\BB_{x,y-p-r+1}\in\DD$
and therefore it should be filled by an element of $\mathcal{A}$
which is larger than $p+r-1$. Consequently, $m_{p+r}(n)>0$.
Now by Lemma \ref{strict}, we have
$$p+r\geq m_1(n)>\cdots>m_{p+r}(n)\geq 1$$
which implies that $m_{k}(n)=p+r+1-k$.
In particular, $m_k(n)>0$ if and only if $k\leq p+r$.
\end{proof}
Let $M=\min\{q+r+1,p+r\}$. By Lemma \ref{semirec2},
for any $1\leq k\leq M$ there are exactly $p+r-k+1$ boxes 
which contain a $k$. 
Lemma \ref{semirec1} determines the position of $p+r-k$ of
these boxes uniquely. Therefore for any $k\leq M$, the location
of precisely one box
containing a $k$ is left to be determined.
In fact it turns out that
after the location of this box is determined, the entire filling is 
also determined uniquely. 
See Lemma \ref{semirec3} below.
\begin{lemma}\label{semirec3}
Let $x'\in\{x,x-1,...,x-p-r+1\}$ and let $y'$ be the smallest
integer such that
$\BB_{x',y'}\in\DD$. Then  

\item[$\bullet$] If $x'\leq x-p-r+M$ then $\BB_{x',y'}$ 
has to be filled with an $x'-x+p+r$.
\item[$\bullet$] For any $y''$ such that 
$y'+1\leq y''\leq y-x'+x-p-r+1$,
the box $\BB_{x',y''}$ has to be filled with a 
$\ov{x'-x+p+r}$.
\end{lemma}
\begin{proof}
Let $b(x')$ be the entry of the box $\BB_{x',y'}$ of the statement
of the lemma. Since the rows and columns are weakly
increasing we have
$$
b(x-p-r+1)\leq b(x-p-r+2)\leq\cdots b(x-p-r+M)
$$
where the inequalities are interpreted in the ordering
of $\mathcal{A}$.
From Definition \ref{amenable}
it follows that $b(x-p-r+j)$ should be unmarked for any
$j\in\{1,...,M\}$. In fact assume on the contrary that 
$b(x-p-r+j)=\ov{l}$ for some $l$. Let $x'=x-p-r+j$
and $y'$ be such that 
$b(x-p-r+j)$ be the entry of the box $\BB_{x',y'}$. 
By Definition \ref{amenable} 
there should exist a box $\BB_{x'',y''}$ such that $y''<y'$
and $\BB_{x'',y''}$ is filled by an $l$. From the shape of 
$\DD$ it follows that $x''>x'$. But then $\BB_{x'',y'}\in\DD$ 
and 
it is impossible to fill the 
latter box such that 
Definition \ref{GSYT} holds.

Since $b(x-p-r+j)$ is unmarked for any $j\in\{1,...,M\}$ and there is exactly
one box filled by a $j$ 
whose location is not determined by Lemma \ref{semirec1},
it turns out that $b(x-p-r+j)=j$ for any such $j$.

Finally, $\BB_{x',y''}$ lies below a box filled by
an $x'-x+p+r-1$ and above a box filled by an $x'-x+p+r$. Consequently,
its entry has to be $\ov{x'-x+p+r}$.

\end{proof}

The proof of Proposition \ref{semirec0} is completed by Lemma
\ref{semirec3}.
The following proposition completes the proof of 
Theorem \ref{main}.

\begin{proposition}\label{llll}
Let $\DD$ be given by $\lambda=\{m>\cdots>1\}$ and
$\mu=\{\mu_1>\cdots>\mu_l\}$ such that $l<m-1$.
Then $\DD$ is a strange diagram.
\end{proposition}

The rest of this section is devoted to the proof of Proposition
\ref{llll}.
Suppose $\DD$ has $p$ columns and $n$ boxes. 
This means that $p=m-l$ and $n={m(m+1)\over 2}-(\mu_1+\cdots+\mu_l)$.
Let $w=w(T)$ be the row word of the filling of $\DD$. $w$ can be expressed as $w=w_1\cdots w_n$ where each $w_i\in\mathcal{A}$.
Since each column of
$\DD$ contains at most one box filled by $1$, $m_1(n)\leq p$.
Lemma \ref{strict} implies that 
\begin{equation}\label{aster}
m_k(n)\leq p-k+1
\textrm{ for any } 
k\in\{1,...,p\}.
\end{equation}
Let $\BB_{x,y}$ denote the box in the lowest row of $\DD$.
(Note that there is only one box in the lowest row.)
For any $i\in\{0,...,p-1\}$ 
let $b_i$ denote the entry of the box $\BB_{x-i,y+i}$. 
Since rows and columns are weakly increasing, the $b_i$'s
form a decreasing sequence in $\mathcal{A}$; i.e.
\begin{equation}\label{eq123}
b_0>\cdots> b_{p-1}.
\end{equation}

\begin{lemma}\label{hard1}
$b_i=p-i$ for any $i\in\{0,...,p-1\}$. 
\end{lemma}
\begin{proof}
First we show that all the $b_i$'s are unmarked.
Suppose on the contrary that for integers 
$i$ and $l$, $b_i=\ov{l}$. Then by Definition 
\ref{amenable}, there should be a box $\BB_{u,v}$
such that $u>x-i$, $v<y-i$ and the entry of $\BB_{u,v}$ is 
$l$. But then it will be impossible to choose an entry for
$\BB_{u,y-i}$ so that Definition \ref{GSYT}
holds. Therefore $b_i$ cannot be marked
for any $i$. 
From equation (\ref{aster}) it is clear that
$b_0\leq p$.

It is now easy to use (\ref{eq123}) to prove that
$b_i=p-i$.
\end{proof}

\begin{definition}\label{diagonal}
A diagonal $D_s$ in $\DD$ is the set of boxes given by
$$D_s=\{\BB_{x',y'}|x'+y'=x+y+s\}\cap \DD.$$
If $D_s\neq\emptyset$ then we can express it as 
$$
D_s=\{\BB_{x_s,y_s},\BB_{x_s+1,y_s-1},...,\BB_{x_s+l_s,y_s-l_s}\}
$$
for some integer $l_s\geq 0$. 
\end{definition}

Our approach is to prove inductively that the entries
of boxes in every $D_s$ are uniquely determined.
The basis of induction is $s=0$, which follows from 
Lemma \ref{hard1}.

\begin{lemma}\label{chalghuz}
Fix $i\in\{0,...,p-1\}$. Let $z>y+i$ be an integer such that
$\BB_{x-i,z}\in\DD$. Then the entry of $\BB_{x-i,z}$ is 
an element of $\mathcal{A}$ strictly less than $p-i$.
\end{lemma}
\begin{proof}
This is because $\BB_{x-i,z}$ lies above $\BB_{x-i,y+i}$ which is
filled by a $p-i$.
\end{proof}

Now let $D_s$ and $D_{s+1}$ be two successive nonempty 
diagonals of 
$\DD$. We can express them as
$$
D_s=\{\BB_{x_s,y_s},...,\BB_{x_s+l_s,y_s-l_s}\}
\textrm{ and } 
D_{s+1}=\{\BB_{x_{s+1},y_{s+1}},...,\BB_{x_{s+1}+l_{s+1},y_{s+1}-l_{s+1}}\}
$$
where $l_s,l_{s+1}\geq 0$.
It can be seen that one of the following cases can happen:

\begin{itemize}
\item[$\bullet$] $x_{s+1}=x_s+1$ and $y_{s+1}=y_s$.
In this case $l_{s+1}=l_s-1$.
\item[$\bullet$] $x_{s+1}=x_s$ and $y_{s+1}=y_s+1$.
In this case $l_{s+1}=l_s$.
\end{itemize}

From Lemma
\ref{hard1} it follows that $b_0=p$, i.e. $m_p(n)$=1. 
Now Lemma \ref{strict} implies that $m_k(n)=p-k+1$. 
The following lemma now follows immediately
from Lemma \ref{chalghuz}. 

\begin{lemma}\label{malghuz}
For any $k$ such that $p\geq k\geq 1$ and 
any $u\geq x-p+k$ there exists 
an integer $v$
such that $\BB_{u,v}$ is filled by a $k$.
\end{lemma}
Next we prove the following lemma, which completes 
the proof of Proposition \ref{llll}.
\begin{lemma}\label{hard2}
The entries of the boxes of $D_{s+1}$ are given
according to the following rules:

\item[$\bullet$]
If $y_{s+1}=y_s+1$ then for any $i\in\{0,...,l_{s+1}\}$ 
the entry of $\BB_{x_{s+1}+i,y_{s+1}-i}$ is $\ov{i+1}$.

\item[$\bullet$]
If $x_{s+1}=x_s+1$ then for any $i\in\{0,...,l_{s+1}\}$ 
the entry of $\BB_{x_{s+1}+i,y_{s+1}-i}$ is $i+1$.

\end{lemma}

\begin{proof}

We use induction on $i$. Namely, we assume that 
the entries of the boxes in $D_s$ are given according to
the rules stated in the lemma, and
then we use induction on $i$
to prove that the entry of $\BB_{x_{s+1}+i,y_{s+1}-i}$
is given according to the rules in the statement of the lemma
as well.

\vspace{1mm}
{\bf Case 1.} 
Obviously $\BB_{x_{s+1},y_{s+1}}=\ov{1}$, since it lies above the box 
$\BB_{x_s,y_s}$
whose entry belongs to $\{1,\ov{1}\}$.
This amounts for the basis of induction, 
i.e. $i=0$. Next suppose $i>0$. Then $\BB_{x_{s+1}+i,y_{s+1}-i}$
lies to the right of $\BB_{x_s+i-1,y_s-i+1}$ and above $\BB_{x_s+i,y_s-i}$.
However, by our assumption about $D_s$, 
the entries of the latter two boxes, which 
belong to $D_s$, are known. A simple argument based on monotonicity of 
rows and columns
implies that if $b$ is the entry of  $\BB_{x_{s+1}+i,y_{s+1}-i}$
then $b\in\{i,\ov{i+1}\}$. It suffices to show that 
$b=i$ is impossible. The proof is by contradiction.
Suppose on the contrary that $b=i$. By induction hypothesis
we know that $\BB_{x_{s+1}+i-1,y_{s+1}-i+1}$ is filled by $\ov{i}$,
which implies that it is impossible to choose an entry for 
$\BB_{x_{s+1}+i,y_{s+1}-i+1}$ such that Definition \ref{GSYT}
holds. Consequently,
$b=\ov{i+1}$.

{\bf Case 2.}
First note that
$\BB_{x_{s+1}+i,y_{s+1}-i}$ lies to the right of $\BB_{x_s+i,y_s-i}$ and 
above $\BB_{x_s+i+1,y_s-i-1}$. But the entries of the latter two boxes are
known, and a simple argument based on monotonicity of rows and columns 
implies that if $b$ is the entry of $\BB_{x_{s+1}+i,y_{s+1}-i}$ then
$b\in\{i+1,\ov{i+2}\}$. Next we show that $b$ cannot be equal to 
$\ov{i+2}$. Suppose on the contrary that $b=\ov{i+2}$. 
Suppose $b$ corresponds to the letter $w_{u_1}$ of the 
row word of the filling. Then by a bookkeeping argument we have
\vspace{-1.8mm}
$$
m_{i+1}(n+u_1)=m_{i+2}(n+u_1).
$$
By Lemma \ref{chalghuz} and Lemma \ref{malghuz}
there exists a $z$ such that the box
$\BB_{x_{s+1}+i,z}$ is filled by an $i+1$. 
Obviously $z>y_{s+1}-i$.
Let the entry of  
$\BB_{x_{s+1}+i,z}$ correspond to $w_u$.
Then there should exist an integer $u'$ such that $u_1<u'<u$ and
$w_{u'}=\ov{i+1}$, because otherwise we will have
$m_{i+1}(n+u-1)=m_{i+2}(n+u-1)$ and $w_u=i+1$ which contradicts
the second property of Definition \ref{amenable}. 

However, existence of $u'$ leads to a contradiction too. In fact 
$w_{u'}$ should correspond to a box $\BB_{x'',y''}$ such that 
either $y''<z$ or $y''=z$ and $x''<x_{s+1}+i$. 
However, if $x''<x_{s+1}+i$ 
and $y''>y_{s+1}-i$
then $\BB_{x'',y''}$ lies
above $\BB_{x_{s+1}+j,y_{s+1}-j}$ for some $j<i$,
and therefore the entry of $\BB_{x'',y''}$
should be strictly less than $\ov{i+1}$.
(Note that if $i>0$ then 
by induction hypothesis $\BB_{x_{s+1}+i-1,y_{s+1}-i+1}=i$
and if $i=0$ then $\BB_{x_{s+1}+i-1,y_{s+1}-i+1}\notin\DD$.)
If $x''=x_{s+1}+i$ and $y''<z$ then $\BB_{x'',y''}$
lies below $\BB_{x_{s+1}+i,z}$ and therefore its
entry should be strictly larger than $i+1$. Finally, if
$x''>x_{s+1}+i$ and $y''<z$, then $\BB_{x'',z}\in\DD$
but it is impossible to choose its entry such that 
Definition \ref{GSYT} holds.
Therefore all of these cases lead to a contradiction.

\end{proof}


\begin{thebibliography}{ZZZZ}

\bibitem[BTW]{wi2} Billera, Louis and Thomas, Hugh and 
van Willigenburg, Stephanie, 
{\it Decomposable compositions, symmetric quasisymmetric functions 
and equality of ribbon Schur functions}, 
preprint (2005), accepted in Adv. Math.


\bibitem[Bes]{wi1} 
Bessenrodt, Christine
On multiplicity-free products of Schur $P$-functions. 
Ann. Comb. 6 (2002), no. 2, 119--124.





\bibitem[St1]{st1}
 Stembridge, John R. {\it Shifted 
tableaux and the projective 
representations of symmetric groups}. 
Adv. Math. 74 (1989), no. 1, 87--134. 

\bibitem[St2]{st2} Stembridge, John R. 
{\it Multiplicity-free products of Schur functions}. 
Ann. Comb. 5 (2001), no. 2, 113--121. 



\bibitem[Wi]{wi3} van Willigenburg, 
Stephanie Equality of Schur and skew Schur functions. Ann. Comb. 9 (2005), no. 3, 355--362.

\end{thebibliography}
\end{document}